\documentclass[11pt]{amsart}

\usepackage{amsmath, amsthm, amssymb, fullpage}

\usepackage{paralist}
\usepackage{enumerate}

\usepackage{hyperref} 

\newcommand{\FS}{\mathrm{FS}}
\newcommand{\an}{\mathrm{an}}
\newcommand{\alg}{\mathrm{alg}}
\newcommand{\cO}{\mathcal{O}}
\newcommand{\bA}{\mathbb{A}}  
\newcommand{\bD}{\mathbb{D}}  
\newcommand{\bP}{\mathbb{P}}  
\newcommand{\bQ}{\mathbb{Q}}  
\newcommand{\bR}{\mathbb{R}}
\newcommand{\bC}{\mathbb{C}} 
\newcommand{\bL}{\mathbb{L}}
\newcommand{\bZ}{\mathbb{Z}}
\newcommand{\sP}{\mathsf{P}}
\newcommand{\sV}{\mathsf{V}}
\DeclareMathOperator{\PGL}{PGL}
\DeclareMathOperator{\tsp}{tsp}
\DeclareMathOperator{\Res}{Res}
\DeclareMathOperator{\Spec}{Spec}
\DeclareMathOperator{\Rat}{Rat}
\DeclareMathOperator{\Sym}{Sym}
\newcommand{\Vv}{\boldsymbol{v}}
\newcommand{\widevec}[1]{\overrightarrow{#1}}
\newcommand{\bN}{\mathbb{N}}
\newcommand{\cU}{\mathcal{U}}
\newcommand{\cH}{\mathcal{H}}
\DeclareMathOperator{\red}{red}
\DeclareMathOperator{\Frac}{Frac}

\theoremstyle{plain}
\newtheorem{theorem}{Theorem}[section]
\newtheorem{lemma}[theorem]{Lemma}
\newtheorem{mainth}{Theorem}
\newtheorem{maincoro}{Corollary}
\theoremstyle{definition}
\newtheorem*{acknowledgement}{Acknowledgement}

\numberwithin{equation}{section}

\makeatletter
\@namedef{subjclassname@2020}{\textup{2020} Mathematics Subject Classification}
\makeatother

\begin{document}

\title[Hybrid limits and degenerating limits]{Hybrid limits and degenerating limits for the canonical equilibrium measures of rational functions}

\author{Y\^usuke Okuyama}
\address{Division of Mathematics, Kyoto Institute of Technology, Sakyo-ku, Kyoto 606-8585 JAPAN}
\email{okuyama@kit.ac.jp}

\keywords{hybrid limit, degenerating limit, canonical equilibrium measure, rational functions, Berkovich projective line}
\subjclass[2020]{Primary 37F10; Secondary 37P50, 12J25}

\date{\today}

\begin{abstract}
For meromorphic families and sequences of rational functions of fixed
degree $>1$ on the complex projective line, we show that, in either setting,
the hybrid limit theorem
and the degenerating limit theorem for canonical equilibrium measures
are equivalent.
\end{abstract}
\maketitle

\section{Introduction}

For every complex rational function $f\in\bC(z)$ of degree $d>1$, the unique maximal entropy measure $\mu_f$
on $\bP^1(\bC)$ exists and coincides with the canonical equilibrium measure 
$\lim_{n\to\infty}d^{-n}(f^n)^*\omega_{\FS}$ of $f$, where $\omega_{\FS}$ is the normalized spherical volume
form on $\bP^1(\bC)$. 
The degeneration of the measures $\mu_f$ as $f$ meromorphically or sequentially degenerates 
in e.g.\ that the zeros and the poles collide
has recently been actively studied, in both the phase space $\bP^1(\bC)$ itself 
and the hybrid fiber space whose fiber is homeomorphic to
either $\bP^1(\bC)$ or the Berkovich projective line $\sP^1_K$ over some reasonable non-archimedean fields $K$. The key point in both formalisms is that there is a non-archimedean rational function $\phi$ of degree $d>1$ on $\bP^1(K)$ associated to $f$ so that the canonical equilibrium measure $\nu_\phi=\lim_{n\to\infty}d^{-n}(\phi^n)^*\delta_{\xi_g}$ of $\phi$ on $\sP^1_K$, where $\xi_g$ is the Gauss point in $\sP^1_K$, 
exists.
Our aim is to contribute to the study of 
such degenerations of dynamics of $f$ by establishing the equivalence of the seemingly independent two degenerating limit theorems on $\mu_f$ 
in those phase space and hybrid space formalisms.

\subsection{Preliminaries}
For any ring $A$ and any integer $e\ge 0$, we set
$V_{e,A}=A^2\otimes_A\Sym^e((A^2)^\vee)$. The homogeneous resultant form
$\Res$ on $V_{e,A}$ defines the open subscheme
\begin{gather*}
 \Rat_{e,A}=D_+(\Res)\subset\bP(V_{e,A})=\bP^{2e+1}_A.
\end{gather*}
For any field $K$, $\Rat_{e,K}$ also denotes the set of all rational functions of degree $e$ over $K$, which is identified with the $K$-points of $D_+(\Res)$. 

For any Banach ring $(A,\|\cdot\|_A)$, let us denote by
\begin{gather*}
 \sP^1_A=(\bP^1_A)^{\an} 
\end{gather*}
the Berkovich projective line over $A$.
The points $x$ of $(\Spec A[z])^{\an}\subset\sP^1_A$ are the multiplicative
seminorms $|\cdot|_x$ on $A[z]$ satisfying $|a|_x\le\|a\|_A$ for every $a\in A$ (for the details, see \cite{Berkovichbook}).
The structure morphism from $\sP^1_A$ to the Berkovich spectrum of $A$ is the analytification of $\bP^1_A\to\Spec A$ and sends $x$ to the seminorm $a\mapsto |a|_x$ on $A$. For any complete non-archimedean field $K$, 
let $k=k_K=K^\circ/K^{\circ\circ}$ be the residue field of $K$, where $K^\circ$ and $K^{\circ\circ}$ denote the valuation ring of $K$ and the maximal ideal of $K^\circ$, respectively.
Let $\sP^1=\sP^1_K$ be the Berkovich projective line over $K$.
The points $\xi\in\sP^1\setminus\{\infty\}$ are the multiplicative seminorms $|\cdot|_\xi$ on $K[z]$ 
extending the absolute value of $K$.
Suppose now that $K$ is also algebraically closed and that the absolute value (or norm) on $K$ is other than the trivial norm
$|\cdot|_0$.
For every $\xi\in\sP^1$, the tangent (or direction) space $T_\xi\sP^1$ is the set of
all components of $\sP^1\setminus\{\xi\}$, and each direction
$\Vv\in T_\xi\sP^1$ is also written as $U(\Vv)$ when it is convenient,
and is also as $\Vv=\widevec{\xi\xi'}$ for any $\xi'\in U(\Vv)$.
When $e>0$, the action of each $\phi\in\Rat_{e,K}$ on $\bP^1=\bP^1(K)$
extends canonically to a finite analytic selfmap of $\sP^1$ of degree $e$.
For every $\xi\in\sP^1$, the Dirac measure $\delta_\xi$ at $\xi$ on $\sP^1$ satisfies
$\phi^*\delta_\xi=\sum_{\xi'\in\phi^{-1}(\xi)}\deg_{\xi'}(\phi)\delta_{\xi'}$,where $\deg_{\xi'}(\phi)$ is the local degree of $\phi$ at $\xi'$. For each $a\in K$ and each $r\in[0,+\infty)$, 
the point $\xi_{a,r}\in\sP^1\setminus\{\infty\}$ 
is the supremum seminorm
\begin{gather*}
 |h|_{\xi_{a,r}}=\sup_{z\in K:|z-a|\le r}|h(z)|,\quad h\in K[z],
\end{gather*}
and is of type I if $r=0$,
of type II if $r\in|K^\times|$, and of type III if
$r\in\bR_{>0}\setminus|K^\times|$.
Every point of type I, II, or III in $\sP^1\setminus\{\infty\}$
is of this form, and the remaining points in
$\sP^1\setminus\{\infty\}$ are of type IV.
In particular, $\xi_{a,0}$ corresponds to the singleton $\{a\}$,
which is identified with $a$.
The point $\infty$ is also regarded as a type I point in $\sP^1$.
The point $\xi_g:=\xi_{0,1}$ is called 
the Gauss point in $\sP^1$. At $\xi_g$, we have the canonical identification
\begin{gather*}
 T_{\xi_g}\sP^1\cong\bP^1(k),\quad
 \widevec{\xi_ga}\leftrightarrow\widehat{a},\quad a\in\bP^1(K),
\end{gather*}
where $\widehat{a}=a\bmod K^{\circ\circ}$ if $a\in K^\circ$, and
$\widehat{a}=\infty$ otherwise. For more details,
see, e.g., the books and the expository article \cite{BR10,BenedettoBook,Jonsson15}.

We equip the field $\bC((t))$ of complex Laurent series 
centered at $t=0$ with the 
$t$-adic norm $|\cdot|_r$, $r\in(0,1)$, 
so that
\begin{gather*}
 \Bigl|\sum_{j=m}^{\infty}a_jt^j\Bigr|_r=r^m
 \quad\text{when }a_m\neq 0;
\end{gather*}
the Levi-Civita field is 
$\bL=\bL_t=\widehat{\overline{\bC((t))}}$.

\subsection{Meromorphic families}
Fix any $r\in(0,1)$. Let $\bC\{t\}\subset\bC[[t]]\subset\bC((t))$ be the ring of convergent power series. Denoting also by $|\cdot|_\infty$ 
the Euclidean norm on $\bC$, we have the Banach ring
$A_r$ of all $\sum_{j\in\bZ}a_jt^j\in\bC((t))$ having the hybrid norm 
$\sum_{j\in\bZ}(\max\{|a_j|_0,|a_j|_\infty\})r^j<\infty$.
In this meromorphic family setting, the {\em hybrid projective line} is 
the Berkovich projective line $\sP^1_{A_r}$ over $A_r$
(\cite[\S\S1.2--1.3]{F});
the closed disk $\overline{\bD}_r:=\{t\in\bC:|t|_\infty\le r\}$ with the analytic topology 
is identified with the Berkovich
spectrum of $A_r$ by assigning to each $t\in\overline{\bD}_r\setminus\{0\}$ the seminorm
$a\mapsto|a(t)|_\infty^{(\log r)/\log|t|_\infty}$ and to $t=0$ the seminorm $a\mapsto |a|_r$ (\cite[Proposition 1.1 and Theorem 1.2]{F}), and we denote by $\pi_r:\sP^1_{A_r}\to\overline{\bD}_r$ 
the structure morphism of $\sP^1_{A_r}$.
For each $t\in\overline{\bD}_r\setminus\{0\}$, Favre's canonical into homeomorphism
$\psi=\psi_r:\bP^1(\bC)\times(\overline{\bD}_r\setminus\{0\})
 \to\sP^1_{A_r}$
in \cite[Theorem 1.2(3)]{F} maps the line $\bP^1(\bC)\times\{t\}$ onto the fiber 
of the structure morphism over $t$
in $\sP^1_{A_r}$, and
the inclusion $A_r\hookrightarrow\bC((t))$ induces the embedding
$\iota_0:\sP^1_{\bC((t))}\hookrightarrow\sP^1_{A_r}$ onto the fiber over $t=0$ 
in $\sP^1_{A_r}$ (\cite[Theorem 1.2(2)]{F}).
Let us also denote by $p_{\xi_g}$
the projection$:\sP^1_{\bL}\to\sP^1_{\bC((t))}$ induced by
the inclusion field morphism $\bC((t))\subset\bL$. 
For each $e\ge 0$, every $\phi\in\Rat_{e,\Frac(\bC\{t\})}$ is regarded not only as 
an element $\phi\in\Rat_{e,\bL}$ but also as a holomorphic family
$(\phi_t)_{0<|t|_\infty\ll 1}$ in $\Rat_{e,\bC}$. We sometimes denote $\phi$ itself by $\phi_t$
when we discuss base changes of $t$.

Our first principal result is the following.

\begin{mainth}\label{th:merlim}
Fix $d>1$ and $r\in(0,1)$. Then 
the following statements are equivalent.
\begin{enumerate}[(i)]
\item\label{item:merH} 
For every $\phi\in\Rat_{d,\Frac(\bC\{t\})}\subset\Rat_{d,\bL}$, we have
$\lim_{t\to 0}\psi_*\mu_{\phi_t}=(\iota_0)_*(p_{\xi_g})_*\nu_\phi$
weakly on $\sP^1_{A_r}$.
\item\label{item:merC} 
For every $\phi$ as in (\ref{item:merH}), if 
$\nu_\phi(\{\xi_g\})=0$, then
\begin{gather*}
 \lim_{t\to 0}\mu_{\phi_t}
 =\sum_{\Vv\in T_{\xi_g}\sP^1_{\bL}}\nu_\phi\bigl(U(\Vv)\bigr)\,\delta_{\Vv}
 \quad\text{weakly on }\bP^1(\bC)\cong T_{\xi_g}\sP^1_{\bL}.
\end{gather*}
\end{enumerate}
\end{mainth}

The statement (\ref{item:merC}) can be made slightly more general.
Setting $\bD=\{t\in\bC:|t|_\infty<1\}$, let $Y$ be a model obtained from $\bP^1(\bC)\times\bD$ by finitely many blow-ups 
at points over $t=0$. Let $\pi_Y:Y\to\bD$ be the composite blow-up
morphism $Y\to\bP^1(\bC)\times\bD$ followed by the second projection, and $Y_t$ denote the fiber of $\pi_Y$ over each $t\in\bD$.
For every $\xi\in\sP^1_{\bC((t))}$, the completed residue field at $\xi$ is
\begin{gather*}
 \cH(\xi):=\widehat{\Frac\bigl(\bC((t))[T]\big/\{P\in\bC((t))[T]:|P|_\xi=0\}\bigr)},
\end{gather*}
where $T$ is the coordinate on an affine chart of
$\bP^1_{\bC((t))}$ whose analytification contains $\xi$ and the
completion is taken with respect to $|\cdot|_\xi$.
The $\cH(\xi)$-point $\Spec\cH(\xi)\to\bP^1_{\bC((t))}$ extends uniquely to a morphism from $\Spec\cH(\xi)^\circ$ to the proper $\bC[[t]]$-model $Y^{\alg}$ obtained by performing on 
$\bP^1_{\bC[[t]]}$ the blow-ups appearing in defining $Y$, 
and we identify $Y_0$ with the set $Y_0^{\alg}(\bC)$ of all $\bC$-points
of the special fiber $Y^{\alg}_0=(Y^{\alg})_0$ of $Y^{\alg}$ as sets.
The reduction map $\red_Y:\sP^1_{\bC((t))}\to Y_0^{\alg}$, which is anti-continuous
(cf.\ \cite[\S4.1, first paragraph]{BPR}),
is defined so that $\xi$ is mapped to the image of the closed
point in $\Spec\cH(\xi)^\circ$ under this extension
of $\Spec\cH(\xi)\to\bP^1_{\bC((t))}$. 
Set
\begin{gather*}
 S(Y):=\red_Y^{-1}\bigl(
 \{\text{generic points of the irreducible components of }Y_0^{\alg}\}\bigr).
\end{gather*}
For $0<|t|_\infty\ll 1$, we also denote by $\mu_{\phi_t}$ the probability measure
$(\bP^1(\bC)\to Y,\ z\mapsto(z,t))_*\mu_{\phi_t}$ on $Y$.

\begin{maincoro}\label{th:mermodel}
The assertions (\ref{item:merH}) and (\ref{item:merC}) of
Theorem \ref{th:merlim} are also equivalent to the following.
\begin{enumerate}[(i)]
\setcounter{enumi}{2}
\item\label{item:merY}
For every $\phi$ as in (\ref{item:merH}) and every model $Y$ as above,
if $((p_{\xi_g})_*\nu_\phi)(S(Y))=0$, then
\begin{gather*}
 \lim_{t\to0}\mu_{\phi_t}
 =\sum_{a\in Y_0^{\alg}(\bC)}
 \bigl((p_{\xi_g})_*\nu_\phi\bigr)\bigl(\red_Y^{-1}(a)\bigr)\,\delta_a
 \quad\text{weakly on }Y.
\end{gather*}
\end{enumerate}
\end{maincoro}

On (\ref{item:merH}), the hybrid limit theorem there was established by Favre
\cite[Theorem B]{F}.
On (\ref{item:merC}), the degenerating limit theorem there
was asserted by DeMarco--Faber \cite[Theorem B]{DF};
the argument in the proof was complemented by \cite[Theorem B]{ODF}.

\subsection{Degenerating sequences}\label{sec:seqsetup}

Recall that the Stone--\v{C}ech compactification $\beta\bN$ of the discrete space $\bN$
is the space of ultrafilters on $\bN$, where each $n\in\bN$ corresponds to
the principal ultrafilter $\{E\subseteq\bN:n\in E\}$ and each $\cU\in\beta\bN\setminus\bN$
is called a non-principal ultrafilter.
For each sequence $\epsilon=(\epsilon_n)_n$ in $(0,1]$ 
tending to $0$ as $n\to\infty$,
we have the Banach ring
$A^\epsilon$ of all $(a_n)_n$ in $\bC^{\bN}$
having the $\epsilon$-norm $\|(a_n)_n\|_\epsilon:=\sup_{n\in\bN}|a_n|_\infty^{\epsilon_n}<\infty$.
For every $\cU\in\beta\bN$, the residual ring
\begin{gather*}
 \cH_\epsilon(\cU):=A^\epsilon\big/\Bigl\{(a_n)_n\in A^\epsilon:
       \lim_{n\to\cU}|a_n|_\infty^{\epsilon_n}=0\Bigr\}
\end{gather*}
is an algebraically closed and complete valued field so that the 
(non-trivial) norm of the class $a_\cU$ represented by $(a_n)_n\in A^\epsilon$
is $|a_\cU|_\epsilon:=\lim_{n\to\cU}|a_n|_\infty^{\epsilon_n}$;
for each $n\in\bN$, $\cH_\epsilon(n)\cong(\bC,|\cdot|_\infty^{\epsilon_n})$ isometrically, and $\cU\in\beta\bN\setminus\bN$ if and only if
$\cH_\epsilon(\cU)$ is non-archimedean. In this sequence setting, the {\em hybrid projective line} is 
the Berkovich projective line $\sP^1_{A^\epsilon}$ over $A^\epsilon$, and
the Stone--\v{C}ech compactification $\beta\bN$ of $\bN$
is identified with the Berkovich spectrum of $A^\epsilon$
by assigning to each $\cU\in\beta\bN$ the seminorm
$(a_n)_n\mapsto|a_\cU|_\epsilon$ (\cite[Theorem 4.8]{FG}).
Let $\pi_\epsilon:\sP^1_{A^\epsilon}\to\beta\bN$ be the structure morphism.
For every $\cU\in\beta\bN$, the quotient map
$A^\epsilon\to\cH_\epsilon(\cU)$ induces the embedding
\begin{gather*}
 \iota_\cU:\sP^1_{\cH_\epsilon(\cU)}\hookrightarrow\sP^1_{A^\epsilon}
\end{gather*}
onto $\pi_\epsilon^{-1}(\cU)$ (\cite[\S\S2.2.1 and 4.2.2]{FG});
for every $n\in\bN$, the quotient map $A^\epsilon\to\cH_\epsilon(n)$
is the evaluation $(a_j)_j\mapsto a_n$ under the identification
$\cH_\epsilon(n)\cong(\bC,|\cdot|_\infty^{\epsilon_n})$,
and on each affine chart of $\sP^1_{\cH_\epsilon(n)}$, $\iota_n$
sends a seminorm $|\cdot|_\xi$ on $\bC[T]$ to the seminorm on $A^\epsilon[T]$ given by
$|\sum_j (a_{j,m})_m T^j|_{\iota_n(\xi)}
 =|\sum_j a_{j,n}T^j|_\xi$, $(a_{j,m})_m\in A^\epsilon$
(\cite[\S2.2.1 and Proposition 4.9(1)]{FG}).
For each $n\in\bN$, we also have the inverse
$s_{\epsilon_n}^{-1}:\bP^1(\bC)\to\sP^1_{(\bC,|\cdot|_\infty^{\epsilon_n})}$
of Favre and Gong's homeomorphism $s_n$ in \cite[\S2.2.2]{FG}, which sends each $z\in\bC$ on each affine chart
to the seminorm $P\mapsto|P(z)|_\infty^{\epsilon_n}$ on $\bC[T]$.
Recall also that for each $\cU\in\beta\bN\setminus\bN$, 
\begin{gather*}
 \Bigl\{\widehat{a_\cU}:a_\cU\in\cH_\epsilon(\cU)
=\bP^1(\cH_\epsilon(\cU))\setminus\{\infty\}
\text{ for some }a=(a_n)_n\in A^\epsilon
\text{ such that }\sup_n|a_n|_\infty<\infty\Bigr\}
\end{gather*}
is a valuation ring having $k_{\cH_\epsilon(\cU)}$ and $\bC$
as its fraction and residue fields, respectively.
The residue homomorphism of this valuation ring sends the mod $\cH_\epsilon(\cU)^{\circ\circ}$ reduction
$\widehat{a_{\cU}}=\widehat{((a_n)_n)_\cU}$ 
to the limit $\lim_{n\to\cU}a_n\in\bC=\bP^1(\bC)\setminus\{\infty\}$ and induces the tangent specialization
\begin{gather*}
 \tsp_{\cU}\colon
 T_{\xi_g}\sP^1_{\cH_\epsilon(\cU)}
 \cong\bP^1\bigl(k_{\cH_\epsilon(\cU)}\bigr)
 \to\bP^1(\bC);
\end{gather*}
see Lemma \ref{th:tsp} and \cite[\S4.2.3, (4.1) and Figure~1]{FG}
for more details.
For every $\cU\in\beta\bN$, each $f\in\Rat_{d,\bZ}(A^\epsilon)$ specializes to $f_{\cU}\in\Rat_{d,\cH_\epsilon(\cU)}$,
which we write $f_n$ when $\cU=n\in\bN$.

Our second principal result is the following sequence version of Theorem \ref{th:merlim}.

\begin{mainth}\label{th:seqlim}
Fix $d>1$. Then the following statements are equivalent.
\begin{enumerate}[(i)]
\item\label{item:seqH} 
For every sequence $\epsilon=(\epsilon_n)_n$ in $(0,1]$
tending to $0$ as $n\to\infty$, 
every $\cU\in\beta\bN\setminus\bN$, and every $f\in\Rat_{d,\bZ}(A^\epsilon)$,
we have $\lim_{n\to\cU}(\iota_n)_*(s_{\epsilon_n})^*\mu_{f_n}=(\iota_\cU)_*\nu_{f_\cU}$
weakly on $\sP^1_{A^\epsilon}$.
\item\label{item:seqC} 
For every $\epsilon$, every $\cU$, and every $f$ as in
(\ref{item:seqH}), if $\nu_{f_{\cU}}(\{\xi_g\})=0$, then
\begin{gather*}
 \lim_{n\to\cU}\mu_{f_n}
 =\sum_{a\in\bP^1(\bC)}\Bigl(\sum_{\Vv\in T_{\xi_g}\sP^1_{\cH_\epsilon(\cU)}:\,\tsp_{\cU}(\Vv)=a}
   \nu_{f_\cU}\bigl(U(\Vv)\bigr)\Bigr)\delta_a
 \quad\text{weakly on }\bP^1(\bC).
\end{gather*}
\end{enumerate}
\end{mainth}

On (\ref{item:seqH}), the sequential hybrid limit theorem there was established by Favre--Gong \cite[Theorem 6.8]{FG}. On (\ref{item:seqC}), the degenerating limit theorem there is an improvement
of the purely atomicness of any weak limit point of $(\mu_{f_n})$
on $\bP^1(\bC)$ \cite[Theorem A]{DF}.
\subsection*{Organization of the paper}

The preliminary lemmas are in Section \ref{sec:prelimlemmas}.
We show Theorem \ref{th:merlim} and Corollary \ref{th:mermodel} in
Section \ref{sec:meromorphic}, and Theorem \ref{th:seqlim} in
Section \ref{sec:sequences}.

\section{Preliminary lemmas}\label{sec:prelimlemmas}

The following, which is a well-known fact, is recorded as a lemma.

\begin{lemma}\label{th:gauss}
Let $K$ be a complete algebraically closed non-archimedean field
with nontrivial absolute value, and let $\phi\in\Rat_{d,K}$ with $d>1$. Then the following statements
are all equivalent: $\nu_\phi(\{\xi_g\})>0$,
$\nu_\phi=\delta_{\xi_g}$, and $\phi^{-1}(\xi_g)=\{\xi_g\}$.
\end{lemma}

\begin{proof}
Noting that $\nu_\phi(\{\phi(\xi_g)\})\ge\nu_\phi(\{\xi_g\})$
by $\phi_*\nu_\phi=\nu_\phi$, 
which follows from $\phi^*\nu_\phi=d\nu_\phi$,
if $\nu_\phi(\{\xi_g\})>0$, then
$\phi^{n+\ell}(\xi_g)=\phi^n(\xi_g)$ for some $n\ge 0$
and some $\ell\ge 1$
since $\nu_\phi$ is a probability measure on $\sP^1$.
Then $d^\ell\cdot\nu_\phi(\{\phi^n(\xi_g)\})
=((\phi^\ell)^*\nu_{\phi})(\{\phi^n(\xi_g)\})
=\deg_{\phi^n(\xi_g)}(\phi^\ell)\nu_\phi(\{\phi^n(\xi_g)\})$
(also by Fubini's theorem) and in turn 
$\phi^{-\ell}(\{\xi_g\})=\{\xi_g\}$, 
which with the defining property $\lim_{n\to\infty}d^{-n}(\phi^n)^*\delta_{\xi_g}=\nu_\phi$ yields $\nu_\phi=\delta_{\xi_g}$. 
The other implications are now immediate. 
\end{proof}

\begin{lemma}\label{th:radial}
Let $K$ be a complete algebraically closed non-archimedean field and pick any dense subset $S$ in $K$. Then the $\bR$-linear span of all the functions
on $\sP^1$, $\xi\mapsto g(|z-b|_\xi)$ for $\xi\ne\infty$ and $\infty\mapsto g(+\infty)$, where $g\in C([0,+\infty])$ and $b\in S$, is uniformly dense in $C(\sP^1)$.
\end{lemma}

\begin{proof}
The $\bR$-linear span of those functions on $\sP^1$
forms an $\bR$-subalgebra in $C(\sP^1)$; indeed, for $a,b\in S$ and $g,h\in C([0,+\infty])$, noting that
$|z-a|_\xi=|z-b|_\xi$ when $\min\{|z-a|_\xi,|z-b|_\xi\}>|a-b|$, we have 
\begin{multline*}
 \bigl(h(|z-b|_\xi)-h(|a-b|)\bigr)\cdot\bigl(g(|z-a|_\xi)-g(|a-b|)\bigr)\\
 =\Bigl(\bigl(h(\cdot)-h(|a-b|)\bigr)\cdot\bigl(g(\cdot)-g(|a-b|)\bigr)\Bigr)
 \bigl(\max\{|z-a|_\xi,|a-b|\}\bigr)\quad\text{on }\sP^1.
\end{multline*}
This span contains all constant functions on $\sP^1$, and moreover 
separates points in $\sP^1$; indeed, for any distinct 
$\xi,\xi'\in\sP^1$, when $\{\xi,\xi'\}\subset\sP^1\setminus\{\infty\}$,
we in turn have $|P|_\xi\neq|P|_{\xi'}$ for some $P\in K[z]$,
$|z-a|_\xi<|z-a|_{\xi'}$ for some $a\in K$ exchanging $\xi,\xi'$ if necessary, and
$|z-b|_\xi<|z-b|_{\xi'}$ for some $b\in S$, so that $1/(1+|z-b|_\cdot)$ separates 
those $\xi,\xi'$, and so does it even when $\infty\in\{\xi,\xi'\}$. 
Hence we are done by the Stone--Weierstrass theorem.
\end{proof}

\begin{lemma}\label{th:tsp}
Let $\epsilon=(\epsilon_n)_n$ be a sequence in $(0,1]$ tending to $0$
as $n\to\infty$ and pick any $\cU\in\beta\bN\setminus\bN$.
Let $\Vv\in T_{\xi_g}\sP^1_{\cH_\epsilon(\cU)}$ be such that
$\tsp_\cU(\Vv)\in\bC$, and let $a^{\Vv}=(a_n^{\Vv})_n\in A^\epsilon$
be such that $\sup_{n\in\bN}|a_n^{\Vv}|_\infty<\infty$ and
$\widehat{(a^{\Vv})_\cU}=\Vv$ under the canonical identification
$\bP^1(k_{\cH_\epsilon(\cU)})
\cong T_{\xi_g}\sP^1_{\cH_\epsilon(\cU)}$.
Then we have $U(\Vv)=\{\xi\in\sP^1_{\cH_\epsilon(\cU)}\setminus\{\infty\}:|z-(a^{\Vv})_\cU|_\xi<1\}$
and $\lim_{n\to\cU}a_n^{\Vv}=\tsp_\cU(\Vv)$.
\end{lemma}

\begin{proof}
Pick $\Vv\in T_{\xi_g}\sP^1_{\cH_\epsilon(\cU)}$ such that
$\tsp_\cU(\Vv)\in\bC$, and fix $a^{\Vv}$ as in the statement.
Then by the first condition on $a^{\Vv}$, we have
$|(a^{\Vv})_\cU|_\epsilon=\lim_{n\to\cU}|a_n^{\Vv}|_\infty^{\epsilon_n}\le 1$,
and in turn $\xi_{(a^{\Vv})_\cU,1}=\xi_{0,1}=\xi_g$,
and by the second condition on $a^{\Vv}$,
$(a^{\Vv})_\cU\in U(\Vv)$. Hence the first assertion holds.
The second one is by the definition of $\tsp_\cU$.
\end{proof}

\section{Proofs of Theorem \ref{th:merlim} and Corollary \ref{th:mermodel}}\label{sec:meromorphic}

\begin{proof}[Proof of Theorem \ref{th:merlim}]
 Let us first assume (\ref{item:merH}),
 and pick any such $\phi$ as in (\ref{item:merH}) and 
 that $\nu_\phi(\{\xi_g\})=0$.
 Pick any weak limit point
 $\mu=\lim_{n\to\infty}\mu_{\phi_{t_n}}$ on $\bP^1(\bC)$
 of $(\mu_{\phi_t})_{0<|t|_\infty\ll 1}$ for some sequence $(t_n)$
 tending to $0$ as $n\to\infty$.
 Let us see that for every $\Vv\in T_{\xi_g}\sP^1_{\bL}$,
 $\nu_\phi\bigl(U(\Vv)\bigr)\le\mu(\{\Vv\})$; 
 for,
 we write $\Vv=\widevec{\xi_g a}$, $a\in\bP^1(\bC)$. 
 When $a\in\bC$, $U(\Vv)=\{\xi\in\sP^1_{\bL}\setminus\{\infty\}:|z-a|_\xi<1\}$.
 For each $\rho\in(0,1)$ and each $m\in\bN$, the positive function
 \begin{gather*}
 x\mapsto \begin{cases}
  \max\bigl\{0,\min\{1,m\cdot (\rho-|z-a|_x)\}\bigr\}
  &\text{if }x\in(\Spec A_r[z])^{\an}\\
  0
  &\text{otherwise}
 \end{cases}\quad\text{on }\sP^1_{A_r}
 \end{gather*}
 is continuous on $\sP^1_{A_r}$ and 
 $\le 1_{\{x\in(\Spec A_r[z])^{\an}:\,|z-a|_x<\rho\}}$.
 Then under the assumption (\ref{item:merH})
 and the choice of $\mu$, for every $\eta\in(0,1)$, 
 noting that $|t_n|_\infty^{(\log\rho)/\log r}<\eta$ 
 for $n\gg 1$ and using the monotone convergence theorem, we have
 \begin{multline}
  \nu_\phi(\{\xi\in\sP^1_{\bL}\setminus\{\infty\}:|z-a|_\xi<\rho\})
  =\lim_{m\to\infty}\int_{\sP^1_{\bL}\setminus\{\infty\}}
   \max\bigl\{0,\min\{1,m\cdot (\rho-|z-a|_\xi)\}\bigr\}\,\nu_\phi\\
  =\lim_{m\to\infty}\lim_{n\to\infty}\int_{\bC}
   \max\bigl\{0,\min\{1,m\cdot (\rho-|z-a|_\infty^{(\log r)/\log|t_n|_\infty})\}\bigr\}\mu_{\phi_{t_n}}\\
  \le\liminf_{n\to\infty}
       \mu_{\phi_{t_n}}\bigl(\bigl\{z\in\bC:|z-a|_\infty^{(\log r)/\log|t_n|_\infty}<\rho\bigr\}\bigr) \le\mu\bigl(\{z\in\bC:|z-a|_\infty\le\eta\}\bigr),
 \label{eq:tube}
 \end{multline}
 and we are done in this case by letting $\rho\nearrow 1$ and then $\eta\searrow 0$.
 When $a=\infty$, $U(\Vv)=\{\xi\in\sP^1_{\bL}\setminus\{\infty\}:|z|_\xi>1\}\sqcup\{\infty\}$, and 
 we are also done by replacing $z-a$ with $1/z$
 in the computation \eqref{eq:tube}. Once this estimate is at our disposal, 
 under the assumption $\nu_\phi(\{\xi_g\})=0$, we have
 $1=\sum_{\Vv\in T_{\xi_g}\sP^1_{\bL}}\nu_\phi(U(\Vv))
 \le\sum_{\Vv\in T_{\xi_g}\sP^1_{\bL}}\mu(\{\Vv\})\le 1$,
 so that $\mu=\sum_{\Vv\in T_{\xi_g}\sP^1_{\bL}}\nu_\phi(U(\Vv))\,\delta_{\Vv}$.
 Hence the degenerating convergence in (\ref{item:merC}) holds.

 Let us next assume (\ref{item:merC}) and 
 pick any such $\phi$ as in (\ref{item:merH}).
 Fix $h\in C(\sP^1_{A_r})$. By the Newton--Puiseux theorem and
 Lemma \ref{th:radial}, for any $\epsilon>0$,
 the pullback $h\circ\iota_0\circ p_{\xi_g}$ of $h$
 has an $\epsilon$-uniform
 approximation by the finite sum of some functions 
 $\xi\mapsto g_j(|z-b_j|_\xi)$ on $\sP^1_{\bL}$, 
 $g_j\in C([0,+\infty])$ and
 $b_j\in\bC[t^{1/m},t^{-1/m}]$, $m\in\bN$, for each $j$.
 Let us see that
 for every pair of
 $g\in C([0,+\infty])$ and $b\in\bC[t^{1/m},t^{-1/m}]$,
 $m\in\bN$,
 \begin{gather}
 \lim_{s\to 0}\int_{\bP^1(\bC)}
 g\bigl(|z-b(s^m)|_\infty^{(\log r)/\log|s^m|_\infty}\bigr)\,\mu_{\phi_{s^m}}(z)
 =\int_{\sP^1_{\bL}}g(|z-b|_\xi)\,\nu_\phi(\xi);
 \label{eq:merint}
 \end{gather}
 indeed, for every $q\in\bQ$ and every $k\in m\bN$ such that $kq\in\bZ$,
 set 
 \begin{gather*}
 M(z):=\frac{1}{t^q}(z-b)\in\PGL\bigl(2,\Frac(\bC\{t^{1/k}\})\bigr)
 \subset\PGL(2,\bL),
 \end{gather*}
 so that $M_{s^k}\circ\phi_{s^k}\circ M_{s^k}^{-1}\in\Rat_{d,\Frac(\bC\{s\})}\subset\Rat_{d,\bL_s}$.

 We note that the base change $t=s^k$ induces
 the homeomorphism $\sP^1_\bL\ni \xi\mapsto y\in\sP^1_{\bL_s}$, where the $s$-adic norm 
 on $\bL_s$ is also $|\cdot|_r$, so that
 when $\xi\in\sP^1_\bL\setminus\{\infty\}$,
 \begin{gather*}
 \biggl|\sum_i a_i(s^k)z^i\biggr|_y
 =\biggl|\sum_i a_i(t)z^i\biggr|_\xi^k,
 \end{gather*}
 so in particular $|z-b(s^k)|_y=|z-b|_\xi^k$, and
 that $y=\infty$ when $\xi=\infty$,
 not only preserving the Gauss points in $\sP^1_{\bL_t}$
 and $\sP^1_{\bL_s}$ but also pushing forward the local degree
 function of $\phi$ to that of $\phi_{s^k}$, 
 so in particular pushing forward the canonical equilibrium measure $\nu_\phi$ on $\sP^1_{\bL}$ of $\phi$ 
 to the canonical equilibrium measure $\nu_{\phi_{s^k}}$ 
 of $\phi_{s^k}$ on $\sP^1_{\bL_s}$. 
 We also note that 
$(M_{s^k})_*\nu_{\phi_{s^k}}=\nu_{M_{s^k}\circ\phi_{s^k}\circ M_{s^k}^{-1}}$ on $\sP^1_{\bL_s}$.

For at most one $q\in\bQ$,
we have $((M_{s^k})_*\nu_{\phi_{s^k}})(\{\xi_g\})>0$
or equivalently, by Lemma \ref{th:gauss},
$(M_{s^k})_*\nu_{\phi_{s^k}}=\delta_{\xi_g}$
on $\sP^1_{\bL_s}$
since ($r\in(0,1)$ and $\xi_g$ is independent of $q$ and) 
$|z-b|_{M^{-1}(\xi_g)}=r^q$ for such 
$q\in\bQ$.
 Unless $q\in\bQ$ is such at most one exception,
 under the assumption (\ref{item:merC}), we compute 
 \begin{multline*}
 \liminf_{s\to 0}\mu_{\phi_{s^m}} \bigl(\bigl\{z\in\bC:|z-b(s^m)|_\infty<|s^m|_\infty^q\bigr\}\bigr)
 =\liminf_{s\to 0}\bigl((M_{s^k})_*\mu_{\phi_{s^k}}\bigr)\bigl(\{z\in\bC:|z|_\infty<1\}\bigr)\\
 \ge \lim_{s\to 0}\int_{\bP^1(\bC)}
 \max\{0,1-|\cdot|_\infty\}
\bigl((M_{s^k})_*\mu_{\phi_{s^k}}\bigr) =\sum_{a\in\bP^1(\bC)}
 \nu_{M_{s^k}\circ\phi_{s^k}\circ M_{s^k}^{-1}}\bigl(U(\widevec{\xi_g a})\bigr)\max\{0,1-|a|_\infty\}\\
 \ge \bigl((M_{s^k})_*\nu_{\phi_{s^k}}\bigr)\bigl(U(\widevec{\xi_g0})\bigr)\cdot 1
+0\\
 =\nu_{\phi_{s^k}}\bigl(\bigl\{y\in\sP^1_{\bL_s}\setminus\{\infty\}:|z-b(s^k)|_y<r^{kq}\bigr\}\bigr)
 =\nu_\phi\bigl(\bigl\{\xi\in\sP^1_{\bL}\setminus\{\infty\}:|z-b|_\xi<r^q\bigr\}\bigr)
 \end{multline*}
and similarly
 \begin{multline*}
 \limsup_{s\to 0}\mu_{\phi_{s^m}}\bigl(\bigl\{z\in\bC:|z-b(s^m)|_\infty<|s^m|_\infty^q\bigr\}\bigr)
 =\limsup_{s\to 0}\bigl((M_{s^k})_*\mu_{\phi_{s^k}}\bigr)\bigl(\{z\in\bC:|z|_\infty<1\}\bigr)\\
 \le \lim_{s\to 0}\int_{\bP^1(\bC)}
 \frac{1}{\max\{1,|\cdot|_\infty\}}\bigl((M_{s^k})_*\mu_{\phi_{s^k}}\bigr)
 =\sum_{a\in\bP^1(\bC)}
 \frac{((M_{s^k})_*\nu_{\phi_{s^k}})(U(\widevec{\xi_g a}))}{\max\{1,|a|_\infty\}}\\
 \le\sum_{a\in\bC}
 \bigl((M_{s^k})_*\nu_{\phi_{s^k}}\bigr)\bigl(U(\widevec{\xi_g a})\bigr)
 \le 1-\bigl((M_{s^k})_*\nu_{\phi_{s^k}}\bigr)\bigl(U(\widevec{\xi_g\infty})\bigr)-0\\
 =\nu_{\phi_{s^k}}\bigl(\bigl\{y\in\sP^1_{\bL_s}\setminus\{\infty\}:|z-b(s^k)|_y\le r^{kq}\bigr\}\bigr)
 =\nu_\phi\bigl(\bigl\{\xi\in\sP^1_{\bL}\setminus\{\infty\}:|z-b|_\xi\le r^q\bigr\}\bigr).
 \end{multline*}
 Hence for every $R>0$, also using the regularity of $\nu_{\phi}$, we have
 \begin{multline*}
 \nu_\phi\bigl(\{\xi\in\sP^1_{\bL}\setminus\{\infty\}:|z-b|_\xi<R\}\bigr)\\
 \le\liminf_{s\to 0}\mu_{\phi_{s^m}}
 \bigl(\bigl\{z\in\bC:|z-b(s^m)|_\infty<|s^m|_\infty^{(\log R)/(\log r)}
 (\Leftrightarrow |z-b(s^m)|_\infty^{(\log r)/\log|s^m|_\infty}<R)\bigr\}\bigr)\\
 \le\limsup_{s\to 0}\mu_{\phi_{s^m}}
 \bigl(\bigl\{z\in\bC:|z-b(s^m)|_\infty<|s^m|_\infty^{(\log R)/(\log r)}
 (\Leftrightarrow |z-b(s^m)|_\infty^{(\log r)/\log|s^m|_\infty}<R)\bigr\}\bigr)\\
 \le\nu_\phi\bigl(\{\xi\in\sP^1_{\bL}\setminus\{\infty\}:|z-b|_\xi\le R\}\bigr),
 \end{multline*}
 so that, since $\nu_\phi$ is a probability measure,
 we have \eqref{eq:merint} for $g=\mathbf{1}_{[0,R)}$
 for all but countably many $R>0$,
 and in turn for every $g\in C([0,+\infty])$.

 The pullback of $h$ under the projection $\sP^1_{A_{r^{1/m}}}\to\sP^1_{A_r}$ induced by $t=s^m$
and the function $y\mapsto\sum_jg_j(|z-b_j(s^m)|_y)$
 are both continuous on $\sP^1_{A_{r^{1/m}}}$, and
 the former is still uniformly $\epsilon$-approximated
 by the latter on $\pi_{r^{1/m}}^{-1}(0)$, 
 and in turn on $\pi_{r^{1/m}}^{-1}(s)$ for $0<|s|_\infty\ll 1$
 by properness of $\pi_{r^{1/m}}$. 
 For $0<|s|_\infty\ll 1$, 
 using also \eqref{eq:merint}, we have
 \begin{align*}
 &\biggl|\int_{\sP^1_{A_r}} h\,\bigl((\psi_r(\cdot,s^m))_*\mu_{\phi_{s^m}}\bigr)
 -\int_{\sP^1_{\bL}} h\circ\iota_0\circ p_{\xi_g}\nu_\phi\biggr|
 =\biggl|\int_{\bP^1(\bC)} h\bigl(\psi_r(z,s^m)\bigr)\mu_{\phi_{s^m}}(z)
 -\int_{\sP^1_{\bL}} h\circ\iota_0\circ p_{\xi_g}\nu_\phi\biggr|\\
 \le&\int_{\bP^1(\bC)} 
 \biggl|h\bigl(\psi_r(z,s^m)\bigr)-\sum_j g_j\bigl(
 |z-b_j(s^m)|_\infty^{(\log r)/\log|s^m|_\infty}\bigr)\biggr|
 \,\mu_{\phi_{s^m}}(z)\\
 &+\sum_j\biggl|\int_{\bP^1(\bC)}g_j\bigl(
 |z-b_j(s^m)|_\infty^{(\log r)/\log|s^m|_\infty}\bigr)\mu_{\phi_{s^m}}(z)-\int_{\sP^1_{\bL}}g_j(|z-b_j|_\xi)\,\nu_\phi(\xi)\biggr|\\
 &+\int_{\sP^1_{\bL}}
 \biggl|(h\circ\iota_0\circ p_{\xi_g})(\xi)-\sum_jg_j(|z-b_j|_\xi)\biggr|
 \nu_\phi(\xi)
 <3\epsilon.
 \end{align*}
 Hence the hybrid convergence in (\ref{item:merH}) holds also noting that $\psi_{r^{1/m}}(z,s)$ is mapped to $\psi_r(z,s^m)$ under the above projection. 
\end{proof}

\begin{proof}[Proof of Corollary \ref{th:mermodel}]

Let us assume (\ref{item:merH}) and pick any $\phi$ and 
any $Y$ such as in (\ref{item:merY}) and that $((p_{\xi_g})_*\nu_\phi)(S(Y))=0$.
Recall also the reduction map $\red_Y$ from $\sP^1_{\bC((t))}$
to the $\bC[[t]]$-model $Y^{\alg}$. Then
$\sP^1_{\bC((t))}\setminus S(Y)
=\bigsqcup_{a\in Y_0^{\alg}(\bC)}\red_Y^{-1}(a)$
(cf.\ \cite[\S1.3, second paragraph]{BPR}),
noting that the special fiber $Y_0^{\alg}$ over $t=0$ is of finite type over $\bC$ and of pure dimension $1$. The models $Y$ and $Y^{\alg}$ are also the analytification over $\bD\subset\bC=(\Spec\bC[t])(\bC)$
and the base change to $\bC[[t]]$, respectively, of the model over $\bC[t]$
obtained by performing on $\bP^1_{\bC[t]}$ the blow-ups defining $Y$.
Similarly, the scheme $\bP^1_{A_r}$ is also the base change of this 
intermediate $\bC[t]$-model to $A_r$;
indeed, the composite blow-up morphism becomes an isomorphism over $A_r$
since the centers of these blow-ups all lie over $t=0$ and $t^{-1}\in A_r$.
We identify the special fiber of this intermediate $\bC[t]$-model
with $(Y^{\alg})_0$, the $\bC$-points of which are identified with $Y_0$
(see Section 1.2).

Let $\tilde{\psi}_{r,Y}\colon\bP^1(\bC)\times(\overline{\bD}_r\setminus\{0\})
\to\pi_Y^{-1}(\overline{\bD}_r\setminus\{0\})$ be the lift of the identity
through the composite blow-up morphism.
Pick any open subset $V$ in $Y$ and set
\begin{gather*}
 \sV:=\sP^1_{A_r}\setminus
 \overline{\psi_r\bigl(\tilde{\psi}_{r,Y}^{-1}(Y\setminus V)\bigr)}, 
\quad\text{so that }\psi_r^{-1}(\sV)\subset\tilde{\psi}_{r,Y}^{-1}(V), 
\end{gather*}
which is open in $\sP^1_{A_r}$. Let us see
\begin{gather*}
\red_Y^{-1}\bigl(V\cap Y_0\bigr)\subset \iota_0^{-1}(\sV);
\end{gather*}
for, pick any $a\in V\cap Y_0=V\cap(Y^{\alg})_0$ and
any $\xi\in\red_Y^{-1}(a)$.
Fix an affine open subset $U$ containing $a$ in the intermediate
$\bC[t]$-model and some generators $h_1,\ldots,h_N$ of $\cO(U)$
as a $\bC[t]$-algebra, all of which vanish at $a$.
The image of the extension to $\Spec\cH(\xi)^\circ$ of
$\Spec\cH(\xi)\to\bP^1_{\bC((t))}$ (see Section 1.2) 
is in the open subset $U_{\bC[[t]]}\subset Y^{\alg}$, so that
$\iota_0(\xi)\in(U_{A_r})^{\an}(\subset\sP^1_{A_r})$
and $\max_j\bigl|(h_j)_{A_r}\bigr|_{\iota_0(\xi)}<1$,
where each $(h_j)_{A_r}=h_j\otimes_{\bC[t]} 1$ 
is the pullback under the projection $U_{A_r}
=U\times_{\Spec\bC[t]}\Spec A_r\to U$ of $h_j$; over $t\in\overline{\bD_r}\setminus\{0\}$, we concretely have $(h_j)_{A_r}=(\psi_r)_*\tilde{\psi}_{r,Y}^*h_j$.
By the continuity of $\pi_r:\sP^1_{A_r}\to\overline{\bD}_r$, we first have
$\lim_{\psi_r(z,t)\to\iota_0(\xi)}t=0$. 
Next, for every $t\in\overline{\bD}_r\setminus\{0\}$,
$\pi_r^{-1}(t)\cap(U_{A_r})^{\an}
=\{\psi_r(z,t):z\in\bP^1(\bC)\text{ and }\tilde{\psi}_{r,Y}(z,t)\in U\}$.
Hence $\tilde{\psi}_{r,Y}(z,t)\in U$ if $\psi_r(z,t)$ is close enough
to $\iota_0(\xi)$ and,
as $\psi_r(z,t)\to\iota_0(\xi)$, in turn
\begin{inparaenum}[(i)]
 \item $\max_j\bigl|(h_j)_{A_r}\bigr|_{\psi_r(z,t)}\to\max_j\bigl|(h_j)_{A_r}\bigr|_{\iota_0(\xi)}<1$, 
 \item $\max\{|t|_\infty,\max_j|h_j(\tilde{\psi}_{r,Y}(z,t))|_\infty\}
 =\max\{|t|_\infty,\max_j\bigl|(h_j)_{A_r}\bigr|_{\psi_r(z,t)}^{(\log|t|_\infty)/\log r}\}\to 0$, 
 \item $\tilde{\psi}_{r,Y}(z,t)\to a\in V$ in $Y$ also noting that
 $(t,h_1,\ldots,h_N):U(\bC)\to\bA^{N+1}(\bC)$ is 
 an into homeomorphism, 
 where $U(\bC)$ is the set of all the $\bC$-points of $U$ and is open
 in the space of all the $\bC$-points of the intermediate $\bC[t]$-model, 
 and finally \item $\iota_0(\xi)\in \sV$.
\end{inparaenum}

Once the above inclusion is at our disposal,
under the assumption (\ref{item:merH}), we have
\begin{multline*}
\liminf_{t\to 0}\bigl((\tilde{\psi}_{r,Y}(\cdot,t))_*\mu_{\phi_t}\bigr)(V)
 \ge \liminf_{t\to 0}\bigl((\psi_r(\cdot,t))_*\mu_{\phi_t}\bigr)(\sV)
 \ge \bigl((\iota_0)_*(p_{\xi_g})_*\nu_\phi\bigr)(\sV)\\
 =\bigl((p_{\xi_g})_*\nu_\phi\bigr)\bigl(\iota_0^{-1}(\sV)\bigr)
 \ge\bigl((p_{\xi_g})_*\nu_\phi\bigr)
 \bigl(\red_Y^{-1}(V\cap Y_0)\bigr)
 =\sum_{a\in V\cap Y_0}
 \bigl((p_{\xi_g})_*\nu_\phi\bigr)\bigl(\red_Y^{-1}(a)\bigr).
\end{multline*}
Hence the degenerating convergence in (\ref{item:merY}) holds
since the both measures in (\ref{item:merY}) are probability ones
also under the assumption $((p_{\xi_g})_*\nu_\phi)(S(Y))=0$.

Finally, the degenerating convergence in (\ref{item:merC}) 
follows from that in (\ref{item:merY}) 
in the case of $Y=\bP^1(\bC)\times\bD$. 
\end{proof}

\section{Proof of Theorem \ref{th:seqlim}}\label{sec:sequences}

\begin{proof}[Proof of Theorem \ref{th:seqlim}]
Let us first assume (\ref{item:seqH}),
and pick any $\epsilon$, $\cU$, and $f$ as in (\ref{item:seqH})
such that $\nu_{f_\cU}(\{\xi_g\})=0$.
Set $\mu:=\lim_{n\to\cU}\mu_{f_n}$ weakly on $\bP^1(\bC)$;
we need neither taking a subsequence nor $\limsup,\liminf$ 
in this sequence setting.
Let us see that for every $a\in\bP^1(\bC)$
and every non-empty finite subset
$E\subset\tsp_{\cU}^{-1}(a)$,
$\sum_{\Vv\in E}\nu_{f_\cU}(U(\Vv))\le\mu(\{a\})$;
for, when $a\in\bC$,
by Lemma \ref{th:tsp}, for every $\Vv\in E$,
$U(\Vv)=\{\xi\in\sP^1_{\cH_\epsilon(\cU)}\setminus\{\infty\}:|z-(a^{\Vv})_\cU|_\xi<1\}$
for some $a^{\Vv}=(a_n^{\Vv})_n\in A^\epsilon$ such that
$\sup_{n\in\bN}|a_n^{\Vv}|_\infty<\infty$ and
$\lim_{n\to\cU}a_n^{\Vv}=a$, and then
for every $\rho\in(0,1)$ and every $m\in\bN$, the positive function
\begin{gather*}
 x\mapsto \begin{cases}
  \max_{\Vv\in E}\max\bigl\{0,\min\{1,m\cdot (\rho-|z-a^{\Vv}|_x)\}\bigr\}
  &\text{if }x\in(\Spec A^\epsilon[z])^{\an}\\
  0
  &\text{otherwise}
 \end{cases}\quad\text{on }\sP^1_{A^\epsilon}
\end{gather*}
is continuous on $\sP^1_{A^\epsilon}$ and
$\le 1_{\bigcup_{\Vv\in E}\{x\in(\Spec A^\epsilon[z])^{\an}:\,|z-a^{\Vv}|_x<\rho\}}$.
Similarly to \eqref{eq:tube}, 
under the assumption (\ref{item:seqH}), for every $\eta>0$,
noting that the $U(\Vv)$, $\Vv\in E$, are pairwise disjoint
and using the monotone convergence theorem, 
we have
\begin{gather*}
 \sum_{\Vv\in E}\nu_{f_\cU}\bigl(\bigl\{\xi\in\sP^1_{\cH_\epsilon(\cU)}
 \setminus\{\infty\}:|z-(a^{\Vv})_\cU|_\xi<\rho\bigr\}\bigr)
 \le\mu\bigl(\{z\in\bC:|z-a|_\infty\le\eta\}\bigr),
\end{gather*}
and we are done in this case
letting $\rho\nearrow1$ and then $\eta\searrow0$.
When $a=\infty$, each $\Vv\in E$ specializes to $0$ in the coordinate $1/z$.
Thus by Lemma \ref{th:tsp}, for each $\Vv\in E$, we write
$U(\Vv)=\{\xi\in\sP^1_{\cH_\epsilon(\cU)}\setminus\{0\}:|1/z-(a^{\Vv})_\cU|_\xi<1\}$
for some $a^{\Vv}=(a_n^{\Vv})_n\in A^\epsilon$ such that
$\sup_{n\in\bN}|a_n^{\Vv}|_\infty<\infty$ and
$\lim_{n\to\cU}a_n^{\Vv}=0$.
For every $\rho\in(0,1)$ and every $m\in\bN$, we use the pullback by
$z\mapsto 1/z$ of the continuous function on $\sP^1_{A^\epsilon}$
defined by the above formula with these $a^{\Vv}$.
Then under the assumption (\ref{item:seqH}), for every $\eta>0$,
noting that $\{n\in\bN:\max_{\Vv\in E}|a_n^{\Vv}|_\infty
+\rho^{1/\epsilon_n}<\eta\}\in\cU$ and using the monotone convergence theorem,
\begin{multline*}
 \sum_{\Vv\in E}\nu_{f_\cU}\bigl(\{\xi\in\sP^1_{\cH_\epsilon(\cU)}
 \setminus\{0\}:|1/z-(a^{\Vv})_\cU|_\xi<\rho\}\bigr)\\
 =\lim_{m\to\infty}\int_{\sP^1_{\cH_\epsilon(\cU)}\setminus\{0\}}
 \max_{\Vv\in E}\max\bigl\{0,\min\{1,m\cdot(\rho-|1/z-(a^{\Vv})_\cU|_\xi)\}\bigr\}
 \,\nu_{f_\cU}\\
 =\lim_{m\to\infty}\lim_{n\to\cU}\int_{\bP^1(\bC)\setminus\{0\}}
 \max_{\Vv\in E}\max\bigl\{0,\min\{1,m\cdot(\rho-|1/z-a_n^{\Vv}|_\infty^{\epsilon_n})\}\bigr\}
 \,\mu_{f_n}\\
 \le\lim_{n\to\cU}\mu_{f_n}\Bigl(\bigcup_{\Vv\in E}
 \{z\in\bP^1(\bC)\setminus\{0\}:|1/z-a_n^{\Vv}|_\infty^{\epsilon_n}<\rho\}\Bigr)
 \le\mu(\{z\in\bP^1(\bC)\setminus\{0\}:|1/z|_\infty\le\eta\}),
\end{multline*}
and we are also done by letting $\rho\nearrow1$ and then $\eta\searrow 0$.
Once the desired inequality is at our disposal, we have $\sum_{\Vv\in\tsp_{\cU}^{-1}(a)}\nu_{f_\cU}(U(\Vv))\le\mu(\{a\})$.
Hence the degenerating convergence in (\ref{item:seqC}) holds
by the same argument as that in the meromorphic family case.

Let us next assume (\ref{item:seqC}), 
and pick any $\epsilon$, $\cU$, and $f$ as in (\ref{item:seqH}).
Pick any $a=(a_n)_n\in A^\epsilon$ and any $R>0$,
set $b:=(R^{1/\epsilon_n})_n\in(A^\epsilon)^\times$,
so that $|b_\cU|_\epsilon=R$, and $M(z):=(z-a)/b\in\PGL(2,A^\epsilon)$, and note that 
$(M_{\cU})_*\nu_{f_{\cU}}=\nu_{M_{\cU}\circ f_{\cU}\circ M_{\cU}^{-1}}$
on $\sP^1_{\cH_\epsilon(\cU)}$.
The points $M_\cU^{-1}(\xi_g)=\xi_{a_\cU,R}$ are distinct for distinct
$R$, and $\nu_{f_\cU}$ has at most countably many atoms. Hence
$\bigl((M_\cU)_*\nu_{f_\cU}\bigr)(\{\xi_g\})=0$ for all but countably
many $R>0$. For every $R>0$ for which this mass is zero, under the
assumption (\ref{item:seqC}),
noting that $\tsp_\cU(\widevec{\xi_g0})=0$ and $\tsp_\cU(\widevec{\xi_g\infty})=\infty$, we have
\begin{multline*}
 \nu_{f_\cU}\bigl(\{\xi\in\sP^1_{\cH_\epsilon(\cU)}\setminus\{\infty\}:|z-a_\cU|_\xi<R\}\bigr)\le\lim_{n\to\cU}\mu_{f_n}\bigl(\{z\in\bC:|z-a_n|_\infty^{\epsilon_n}<R\}\bigr)\\
\le\nu_{f_\cU}\bigl(\{\xi\in\sP^1_{\cH_\epsilon(\cU)}\setminus\{\infty\}:|z-a_\cU|_\xi\le R\}\bigr).
\end{multline*}
Now the hybrid convergence in (\ref{item:seqH}) holds 
by an argument similar to and slightly simpler than that in the meromorphic family case.
\end{proof}

\begin{acknowledgement}
 This research was partially supported by JSPS Grant-in-Aid for Scientific Research (C), 23K03129, and
by the Research Institute for Mathematical Sciences,
an International Joint Usage/Research Center at Kyoto University.
\end{acknowledgement}

\end{document}